\documentclass[12pt]{article}
\usepackage{amsmath}
\usepackage{amssymb}
\usepackage{amscd}
\usepackage{amsthm}
\usepackage{amsfonts}
\usepackage{graphicx}
\usepackage{makeidx}         
\usepackage{tabularx} 
\usepackage{hhline}
\usepackage{vmargin} 

\setpapersize[portrait]{A4}
\setmarginsrb{25mm} 
{30mm} 
{20mm} 
{30mm} 
{0mm} 
{0mm} 
{20mm} 
{10mm} 
\newcommand{\CO}{\mathcal{O}}

\newcommand{\BC}{\mathbb C}
\newcommand{\BF}{\mathbb F}

\newcommand{\BN}{\mathbb N}

\newcommand{\BQ}{\mathbb Q}
\newcommand{\BR}{\mathbb R}

\newcommand{\BZ}{\mathbb Z}

\newcommand{\fo}{\mathfrak o}

\newcommand{\ol}{\overline}

\newcommand{\diag}{\textrm{Diag }}

   {\setlength{\parskip}{0pt}%
   \begin{enumerate}%
    }%
   {\end{enumerate}}

\newtheorem{Satz2}{General Criterion2}
\newtheorem{Satz}{General Criterion}
\newtheorem{Thm}[Satz]{Theorem}

\newtheorem{Remark}[Satz]{Remark}

\newtheorem{Fact}[Satz2]{Fact}

\begin{document}
\begin{center}
\textbf{\Large{Definable henselian valuation rings}} \\
\strut \\
Alexander Prestel\label{F1} \\
\end{center}
\strut \\
\strut \\

\textbf{Abstract} We give model theoretic criteria for $\exists
\forall$ and $\forall \exists$- formulas in the ring language to
define uniformly the valuation rings $\CO$ of models $(K, \CO)$ of
an elementary theory $\Sigma$ of henselian valued fields. As one of
the applications we obtain the existence of an $\exists
\forall$-formula defining uniformly the valuation rings $\CO$ of
valued henselian fields $(K, \CO)$ whose residue class field $k$ is
finite, pseudo-finite, or hilbertian. We also obtain $\forall
\exists$-formulas $\varphi_2$ and $\varphi_4$ such that $\varphi_2$
defines uniformly $k[[t]]$ in $k((t))$ whenever $k$ is finite or the
function field of a real or complex curve, and $\varphi_4$
does the job if $k$ is any number field. \\

\section{Introduction}

Let $(K, \CO)$ be a field $K$ together with a valuation ring $\CO$.
We call $\CO$ \textit{definable} in the ring language $L$ if there
exists an $L$-formula $\varphi (x)$ with $x$ a only free variable
and no parameter from $K$ such that
$$
\CO = \{ a \in K \ | \ \varphi (a) \textrm{ hold in } K \}.
$$
We shall mainly be interested in $L$-formulas $\varphi (x)$ of the
following types: \\
$$
\begin{array}{rcl} \exists\textrm{- formula} & : & \exists y_1
\cdots y_n \chi (x, \ol{y} ) \\
\forall\textrm{- formula} & : & \forall y_1 \cdots y_n \chi (x,
\ol{y} ) \\
\exists \forall\textrm{- formula} & : & \exists y_1 \cdots y_n
\forall z_1 \cdots z_m \chi (x, \ol{y}, \ol{z}) \\
\forall \exists\textrm{- formula} & : & \forall y_1 \cdots y_n
\exists z_1 \cdots z_m \chi (x, \ol{y}, \ol{z})\end{array}.
$$
The whole investigation generalizes straight forward to complexer
quantifier types (and is left to the interested reader). We shall
prove (in Section 2) and apply (in Section 3 and 4) the following
model theoretic criteria. \\

\textbf{Characterization Theorem} \textit{Let $\Sigma$ be a first
order axiom system in the ring language $L$ together with a unary
predicate $\CO$. Then there exists a $L$-formula $\varphi (x)$,
defining uniformly in every model $(K, \CO)$ of $\Sigma$ the set
$\CO$, of quantifier type}
\begin{equation*}
\left.
\begin{array}{rcl} \exists & \textrm{ iff } & (K_1 \subseteq K_2
\Rightarrow \CO_1 \subseteq \CO_2) \\
\forall & \textrm{ iff } & (K_1 \subseteq K_2 \Rightarrow \CO_2 \cap
K_1
\subseteq \CO_1) \\
\exists \forall & \textrm{ iff } & (K_1
\stackrel{\exists}{\subseteq} K_2
\Rightarrow \CO_1 \subseteq \CO_2) \\
\forall \exists & \textrm{ iff } & (K_1
\stackrel{\exists}{\subseteq} K_2 \Rightarrow \CO_2 \cap K_1
\subseteq \CO_1) \end{array} \right\} \textrm{ for all models }
(K_1, \CO_1), (K_2, \CO_2) \textrm{ of } \Sigma.
\end{equation*}

\textit{Here $K_1 \stackrel{\exists}{\subseteq} K_2$ means that
$K_1$ is existentially closed in $K_2$, i.e. every $\exists$-formula
$\varrho (x_1, \ldots, x_m )$ with parameters from $K_1$ that holds
in $K_2$
also holds in $K_1$. }\\

If $K_1$ and $K_2$ are fields, this implies that $K_1$ is relatively
algebraically closed in $K_2$. In particular, if $(K_2, \CO_2 )$ is
henselian, then also $(K_1 , \CO_2 \cap K_1)$ is henselian in $K_1$.
Thus if $(K_1, \CO_1)$ is also henselian, we can apply the theory of
henselian valuation rings on a field $K_1$ as explained in Section
4.4 \cite{E-P}. This will yield a series of applications in Sections
3 and 4. \\

It is important to realize that the model theoretic criteria above
do not give explicit $L$-formulas, rather only their existence. But
the knowledge of the existence may help to construct such a formula.
In many cases explicit formulas are already known. Let us mention
here the papers \cite{C-D-L-M}, \cite{A-K}, and \cite{F}. These
papers actually
inspired us to look for general model theoretic criteria. \\

The author is grateful to A. Fehm and J. Schmid for helpful
discussions. \\

\section{Proof of the Characterization Theorem}

Let $\Sigma$ be a first order axiom system in the ring language $L$
enlarged by a unary predicate $\CO$. Moreover, fix a constant $c$.
We denote by $L( \CO, c)$ the enlarged language. An $L( \CO,
c)$-structure then looks like $(K, \CO, a )$ where $K$ is an
$L$-structure, $\CO \subseteq K$ and $a \in K$. Next let $\Phi$ be a
subset of $L (c)$-sentences, i.e. formulas in the ring language $L$
enlarged by $c$ without free variables. We assume that $\Phi$ is
closed by $\wedge$ and $\vee$. Examples of interest to us are the
sets of $\exists, \forall, \exists \forall$, and $\forall
\exists$-sentences in $L(c)$. We furthermore use the following
abbreviation for subsets $\Gamma$ of $L (\CO, c)$-sentences: if
$(K_1 , \CO_1 , a_1 )$ and $(K_2 , \CO_2 , a_2 )$ are two $L (\CO,
c)$-structures and every $\gamma \in \Gamma$ that holds in $(K_1 ,
\CO_1 , a_1 )$ also holds in $(K_2 , \CO_2 , a_2 )$ we write
$$
(K_1 , \CO_1 , a_1 ) \stackrel{\Gamma}{\rightsquigarrow} (K_2 ,
\CO_2 , a_2 ).
$$
Now Lemma 3.1.6 of \cite{P-D}, an easy consequence of the
Compactness Theorem for 1-order logic, immediately gives:

\begin{Thm} Assume that for all models $(K_i ,
\CO_i , a_i )$ of $\Sigma (i = 1,2)$ we have the following
implications:
\newline If $(K_1 , \CO_1 , a_1 )
\stackrel{\Phi}{\rightsquigarrow} (K_2 , \CO_2 , a_2 )$ then $(K_1 ,
\CO_1 , a_1 ) \stackrel{\{
 c \in \CO \}}{\rightsquigarrow} (K_2 , \CO_2 , a_2 )$. Then there
 exists some $\varphi (c) \in \Phi$ such that
 $$
\forall x (x \in \CO \Leftrightarrow \varphi (x))
 $$
 holds in all models $(K, \CO)$ of $\Sigma$, i.e., $\varphi$ defines
 $\CO$ in $K$.
\end{Thm}

The philosophy behind Lemma 3.1.6 is: if a sentence behaves like all
$\varphi \in \Phi$, it is equivalent to some fixed $\varphi \in
\Phi \mod \Sigma$. \\

\textit{Proof of the Characterization Theorem}: Let us first observe
that the right hand side of the equivalences follow clearly from the
corresponding definabilities. \\

\textbf{$\exists$-case:} Assume that for all models $(K_i , \CO_i)
(i = 1,2)$ of $\Sigma$ we have: $K_1 \subseteq K_2 \Rightarrow \CO_1
\subseteq
\CO_2$. \\

Considering two models $(K_1 , \CO_1 , a_1 )$ and $(K_2 , \CO_2 ,
a_2 )$ of $\Sigma$ satisfying
$$
(K_1 , \CO_1 , a_1 ) \stackrel{\Phi}{\rightsquigarrow} (K_2 , \CO_2
, a_2 )
$$
we have to show if $a_1 \in \CO_1$ then also $a_2 \in \CO_2$. We
assume that $a_2 \not = \CO_2$ and consider the set
$$
\Pi = \diag (K_1 , a_1, (b)_{b \in K_1}) \cup Th (K_2, \CO_2, a_2)
$$
of sentences. We claim that $\Pi$ is consistent. In fact, if $\Pi$
would be inconsistent, there would exist some elements $b_1, \ldots,
b_n \in K_1$ and a quantifier free $L$-formula $\chi$ such that
$\chi (a_1, b_1, \ldots, b_n) \in \diag$ of $K_1$ and $\{ \chi \}
\cup Th (K_2, \CO_2, a_2)$ is inconsistent. Hence the
$\exists$-formula $\varphi (c) \equiv \exists b_1, \ldots, b_n \chi
(c, \ol{b})$ of $L(c)$ \footnote{Now $b_1 \cdots b_n$ play the role
of variables.} holds in $(K_1, a_1)$ and we have
$$
Th (K_2 , \CO_2 , a_2 ) \vdash \forall b_1 \cdots b_n \neg \chi (c,
\ol{b}).
$$
This is impossible as $\varphi (c)$ carries over from $(K_1, a_1)$
to $(K_2, a_2)$. \\

Therefore $\Pi$ is consistent and hence has a model
$$
(K^*_2, \CO^*_2, a^*_2)
$$
that is elmentarily equivalent to $(K_2 , \CO_2 , a_2 )$ and
contains an isomorphic copy of $(K_1, a_1)$. After identifying
$(K_1, a_1)$ with its image in $(K^*_2 , \CO^*_2 , a^*_2 )$ we
obtain $K_1 \subseteq K^*_2, a_1 \in \CO_1$, and $a_1 = a^*_2 \not
\in \CO^*_2$. Hence $\CO_1 \not \subseteq \CO^*_2$. This contradicts
our assumption, as
$(K^*_2, \CO^*_2)$ is also a model of $\Sigma$. \\

\textbf{$\exists \forall$-case:} Assume that for all models $(K_i,
\CO_i)$ of $\Sigma$ we have: $K_1 \stackrel{\exists}{\subseteq} K_2
\Rightarrow \CO_2 \subseteq \CO_2$. Looking at the proof of the
$\exists$-case, the only change we need is that $(K_1, a_1)$ (after
identification with its image) is existentially closed in $(K^*_2,
a^*_2)$. this is obtained replacing $\Pi$ by the set
$$
\Pi_{\forall} = Th_{\forall} (K_1, a_1, (b)_{b \in K_1}) \cup Th
(K_2, \CO_2, a_2)
$$
where $Th_{\forall} (K_1, a_1 (b)_{b \in K_1})$ consists of all
$\forall$-formulas
$$
\varphi (x, \ol{b}) \equiv \forall y_1 \cdots y_n \chi (c, \ol{y},
\ol{b}),
$$
where $\chi$ is quantifier free, that hold in $(K_1, a_1, (b)_{b \in
K_1})$. \\

\textbf{$\forall$-case and $\forall \exists$-case:} Is obtained from
the $\exists$-case and the $\exists \forall$-case just by replacing
the sets $\Pi$ and $\Pi_{\forall}$ by
$$
\Pi' = \diag (K_2, a_2, (b)_{b \in K_2}) \cup Th (K_1, \CO_1, a_1)
$$
and
$$
\Pi'_{\forall} = Th_{\forall} (K_2, a_2, (b)_{b \in K_2}) \cup Th
(K_1, \CO_1, a_1)
$$
respectively. \hfill{$\square$} \\

\section{$\exists \forall$-definable henselian valuation rings}

In our applications we shall concentrate here on henselian valued
fields $(K, \CO)$. The maximal ideal of $\CO$ is denoted by $M$, the
residue class field (r.c.f.) by $k = \CO \diagup M$, and the value
group by $v(K)$. If we deal with several valued fields $(K_i,
\CO_i)$ we use corresponding indices for the r.c.f. $k_i$ and the
value groups $v_i (K_i)$. Of particular interest will be the
henselian valuation ring $k [[t]]$ of the fields $k((t))$ of formal
Laurent series and $p$-adic number fields.

Before we proceed to concrete results let us quote some facts about
henselian valued fields from \cite{E-P}. \\

A valued field $(K, \CO)$ is called \textit{henselian} if the
valuation ring $\CO$ of $K$ extends uniquely to the separable
closure $K^s$ of $K$. Note that the trivial valuation $\CO = K$
always is henselian, its residue class field is $K$. Two valuation
rings $\CO_1$ and $\CO_2$ of the same field $K$ are called
comparable if $\CO_1 \subseteq \CO_2$ or $\CO_2 \subseteq \CO_1$;
the upper
one is called \textit{coarser}. Here are some important facts. \\

\begin{Fact} If $\CO_1$ and $\CO_2$ are henselian on $K$ and at
least one of the r.c.f. is not separably closed then $\CO_1$ and
$\CO_2$ are comparable. If $\CO_1$ and $\CO_2$ are not comparable,
then the r.c.f. of $\CO_1, \CO_2$, and of the smallest common
coarsening of $\CO_1$ and $\CO_2$, all are separably closed.
\end{Fact}
\vspace{-6pt}(Theorem 4.4.2 in \cite{E-P}). \\

Now let $\CO_1$ and $\CO_2$ be comparable valuation rings of $K$,
say $\CO_1 \subseteq \CO_2$. (Hence $M_2 \subseteq M_1$.) Then $\fo
= \CO_1 \diagup M_2$ is a valuation ring of the residue class field
$k = \CO_2 \diagup M_2$. Then we obtain from Section 2.3 in
\cite{E-P}
and Corollary 4.1.4: \\

\begin{Fact} The value group of $(k, \fo)$ is isomorphic to a convex
subgroup $\Delta$ of $v (K_1)$ and $v_2 (K_2) \cong v_1 (K_1)
\diagup \Delta$.
\end{Fact}

\begin{Fact} $(K_1, \CO_1)$ is henselian if and only if $(K_1,
\CO_2)$ and $(k, \fo)$ are both henselian.
\end{Fact}

Now let us consider a first order axiom system $\Sigma$ for
henselian valued fields $(K, \CO)$ such that the r.c.f $k= \CO
\diagup M$
\begin{itemize} \item[(1)] is not separably closed \item[(2)] does
not carry a proper henselian valuation.
\end{itemize}
We shall then prove for any two models $(K_1, \CO_1)$ and $(K_2,
\CO_2)$ of $\Sigma$:
$$
K_1 \stackrel{\exists}{\subseteq} K_2 \Rightarrow \CO_1 \subseteq
\CO_2.
$$
Then by the Characterization Theorem there exists an $\exists
\forall$-formula $\varphi (x)$ in the ring language that defines
$\CO$ in every model $(K, \CO)$ of $\Sigma$.

For the proof assume that $K_1 \stackrel{\exists}{\subseteq} K_2$
and $\CO_1 \nsubseteq \CO_2$. As $K_1$ is separably closed in $K_2$,
it follows that $\CO := K_1 \cap \CO_2$ is a henselian valuation
ring of $K_1$. Since by (1) the r.c.f. of $\CO_1$ is not separably
closed, Fact 1 implies $\CO \subsetneqq \CO_1$. Now Fact 3 implies
that $\CO \diagup M_1$ is a proper henselian valuation of $\CO_1
\diagup M_1$. This contradicts (2). Hence $\CO_1 \subseteq \CO_2$,
and we are done.

As an application we obtain
\setcounter{Satz}{0} \begin{Thm} There
is an $\exists \forall$-formula $\varphi (x)$ defining uniformly the
valuation rings of henselian fields $(K, \CO)$ if the residue class
field $k$ of $\CO$ is finite, pseudo-finite, or hilbertian.
\end{Thm}

\textit{Proof:} the class of finite and pseudo-finite fields is the
model class of the theory of finite fields, hence an elementary
class. The class of hilbertian fields is as well elementary. The
union of the elementary classes is again elementary. Let $\Sigma'$
be a first order axiom system for the union. Then let $\Sigma$
express the fact that its models $(K, \CO)$ are henselian valued
fields (not excluding the trivial valuation) such that the r.c.f.
$\CO \diagup M$ satisfies the axioms for $\Sigma'$. \\

We have to check (1) and (2) from above. (1) is clear. (2) is
clearly true for finite fields $k$. If $k$ is pseudo-finite it is a
PAC-field (see \cite{A2}, Section 6, Lemma 2), and PAC-fields do not
carry a proper henselian valuation (unless they a separably closed).
This old result of the author can be found in \cite{F-J}, Corollary
11.5.5. Thus it remains to prove that a hilbertian field $k$ does
not allow a proper henselian valuation ring $\fo$. For contradiction
assume $\fo \subsetneqq k$ is a henselian valuation ring of $k$. We
then choose a separable polynomial $f(x) \in k [X]$ without zero in
$k$. As $\fo$ is henselian, the set $f(k)$ stays away from $0$, say
$f(k) \cap m = \Phi$ ($m$ the maximal ideal of $\fo$). We then
choose $\pi \in m \smallsetminus \{ 0 \}$ and consider the
polynomial
$$
g (X, Y) = f (X) Y^2 + f(X) Y + \pi.
$$

Replacing $Y$ by $Z^{-1}$ and applying Eisenstein, we see that $g
(X,Y)$ is absolutely irreducible. Now let $x$ be any element of $k$.
Then $Y^2 + Y + \frac{\pi}{f(x)}$ maps to $Y(Y+1)$ in $\CO \diagup
M$. Now by Hensel's Lemma $g(x,Y)$ has a zero in $k$, thus is not
irreducible. This contradicts the assumption that $k$ is hilbertian.
\hfill{$\square$} \\

Theorem 1 covers all completions of finite number fields, i.e.
finite extension of the $p$-adic number fields $\BQ_p$ for any prime
$p$. Moreover, it covers all fields $k_0 ((t))$ of Laurent series
with $k_0$ finite, pseudo-finite, or hilbertian. It even covers any
such field $k_0$ together with the trivial valuation. Thus the
sentence $\forall x
\varphi (x)$ is true in all such fields (!). \\

\setcounter{Satz}{0} \begin{Remark} In \cite{A1} Ax gives a $\exists
\forall \exists \forall$-formula that defines $k [[t]]$ uniformly in
$k((t))$ for all fields $k$. Theorem 1 gives an improvement in case
$k$ is finite, pseudo-finite, or hilbertian. In the next Section we
shall consider classes of fields $k$ for which $k[[t]]$ is uniformly
$\forall \exists$-definable. We shall also explain an example $k^*$
of A. Fehm for which $k^* [[t]]$ is not $\forall\exists$-definable
in $k^* ((t))$.
\end{Remark}

\section{$\forall \exists$-definable henselian valuation rings}

For our next theorem we shall need some preparation. As usual we
consider henselian valued fields $(K, \CO)$. This time, however, we
shall require that the value group $v(K)$ does not admit a convex
2-divisible subgroup $\Delta \not = \{ 0 \}$, like discrete value
groups do. This property is easily expressed in the elementary ring
language $L$ together with a predicate $\CO$. In fact, $v(K)$ is
order-isomorphic to $K^{\times} \diagup \CO^{\times}$, and
expressing the existence of a proper convex, 2-divisible subgroup of
$(v(K), \leq)$ can be done by saying in $(v (K), \leq)$:
$$
\exists \gamma (0 < \gamma \wedge \forall \delta (0 \leq \delta \leq
\gamma \Rightarrow \exists \varepsilon \ \delta = 2 \varepsilon )).
$$
Next we need to talk about the $u$-invariant of the residue field
$k$ of $(K, \CO)$ and again be able to do this is first order logic.
We shall make use of the language of quadratic forms as found e.g.
in \cite{E-P}, Section 6.3 or in \cite{P-D}, Chapter 3.

The $u$-\textit{invariant} of a field $k$ is defined to be the
maximal dimension $n \in \BN \cup \{ 0 \}$ of an anisotropic,
quadratic form $\varrho = < a_1, \ldots, a_n >$ with $a_i \in k
\smallsetminus \{ 0 \}$ that has total signature zero (see
\cite{E-K-M}, Chapter VI). For example, $u(k) = 4$ for every finite
number field, and $u (\BC) = 1$. The $u$-invariant of a real field
clearly has to be even, e.g. $u (\BR) = 0$. Here $\varrho$ is called
to be of \textit{total signature} zero over $k$, if for any ordering
$\leq$ of $k$, one half of the $a_i$ is positive and the other half
is negative. There is a quantifier free formula $\zeta (a_1, \ldots,
a_n)$ expressing in the real closure of $(k, \leq)$ that $\varrho$
is not of signature zero. Thus we have to say that there does not
exist an ordering $\leq$ of $k$ such that $\zeta$ holds in $(k,
\leq)$. Using the theory of pre-orderings (see \cite{P-D}) this can
be done in the language of $k$ in case the Pythagoras number $P (k)$
is finite. The Pythagoras number $P (k)$ is the smallest $m \in \BN
\cup \{ 0 \}$ such that every sum of squares in $k$ equals a sum of
$m$ squares. Clearly, if $u(k)$ is finite, then also $P(k)$ is
finite. \\

\begin{Thm} For every non-zero $n \in \BN$ there is an $\forall \exists$-formula $\varphi_n (x)$
defining uniformly the valuation ring $\CO$ of henselian fields $(K,
\CO)$ if the value group $v (k)$ does not admit a convex 2-divisible
subgroup $\Delta \not = \{ 0 \}$, $\textrm{char } \CO \diagup M \not
= 2$, and the $u$-invariant of the residue class field $\CO \diagup
M$ is $n$.
\end{Thm}

\textit{Proof:} We want to apply the Characterization Theorem to the
models $(K, \CO)$ of a first order axiom system $\Sigma$ expressing
that the value group $v (k)$ does not admit a convex 2-divisible
subgroup $\Delta \not = \{ 0 \}$ and that $u (\CO \diagup M) = n$.
Let $(K_1, \CO_1)$ and $(K_2, \CO_2)$ be models of $\Sigma$ and let
$K_1 \stackrel{\exists}{\subseteq} K_2$. We then have to prove that
$\CO
:= K_1 \cap \CO_2 \subseteq \CO_1$. \\

Since $K_1$ is existentially closed in $K_2$ it follows that $\CO$
is henselian and thus $\CO$ and $\CO_1$ are comparable. If not, Fact
1 implies that the r.c.f. of $\CO_1, \CO,$ and the smallest common
coarsening $\CO'$ of $\CO_1$ and $\CO$, all have separably closed
r.c.f. As $\CO_1 \subsetneqq \CO' , \CO_1 \diagup M'$ is a proper
henselian valuation ring of $\CO' \diagup M'$, that has a value
group $\Delta$, not divisible by 2. This follows from the assumption
of the theorem and Fact 2. On the other hand as $\CO' \diagup M'$ is
separably closed, the value group $\Delta$ of $\CO_1 \diagup M'$ has
to be divisible, a contradiction. Therefore $\CO$ and $\CO_1$ are
comparable. It thus remains to exclude $\CO_1 \subsetneqq \CO$.
\\

Let us assume $\CO_1 \subsetneqq \CO$. Then $\CO_1 \diagup M = \fo$
is a proper henselian valuation ring on $\CO \diagup M$ by Fact 3.
The value group $\Delta$ of $\fo$ is a convex subgroup $\not = \{ 0
\}$ of $v_1(K_1)$ and hence by $\Sigma$ not 2-divisible. The r.c.f.
of $\fo$ equals that of $\CO_1$ (Fact 2). As $u (\CO_1 \diagup M_1)
= n$ there exists a quadratic form $\ol{\varrho} = < a_1 + M,
\ldots, a_n + M >$ with $a_i \in \CO^{\times}$ such that
$\ol{\varrho}$ is of total signature zero but not isotropic in
$\CO_1 \diagup M_1$. We then choose some $b \in \CO^{\times}$ such
that with respect to the valuation of $\fo$ its value is not
2-divisible in $\Delta$. Then one can easily check that the
quadratic form
$$
\ol{\varrho}_b := < a_1 + M, \ldots, a_n + M, a_1 b + M, \ldots, a_n
b + M >
$$
cannot be isotropic in $\CO \diagup M$. Moreover, $\ol{\varrho}_b$
is of total signature zero in $\CO \diagup M$. At this point we use
Lemma 4.3.6 and Theorem 2.2.5 of \cite{E-P} to see that every
ordering of $\CO \diagup M$ maps to some ordering of the r.c.f.
$\CO_1 \diagup
M_1$ of $\CO_1 \diagup M$. \\

On the other hand, by $u (\CO_2 \diagup M_2) = n$ we know that
$\ol{\varrho}_b$ is isotropic in the extension $\CO_2 \diagup M_2$
of $\CO \diagup M$. As $(K_2, \CO_2)$ is henselian and $\textrm{char
} \CO_2 \diagup M_2 \not = 2$, it follows from Hensel's Lemma that
the quadratic form
$$
\varrho_b := < a_1, \ldots, a_n, a_1 b, \ldots, a_n b >
$$
is isotropic in $K_2$. Now, since $K_1$ is existentially closed in
$K_2, \varrho_b$ is also isotropic in $K_1$. This, however, clearly
implies that $\ol{\varrho}_b$ is isotropic in $\CO \diagup M$. This
contradiction implies that $\CO_1 \subsetneqq \CO$ cannot hold.
Hence $\CO_2 \cap K_1 = \CO \subseteq \CO_1$, and we are done. \hfill{$\square$} \\

As an application of Theorem 2 we see that the henselian valuation
rings $k[[t]]$ are uniformly definable in $k ((t))$ by some $\forall
\exists$-formula $\varphi_n$ in case\footnote{As $u (\BR)= 0$, this
theorem does not cover the case of $\BR[[t]]$. Replacing, however,
the $u$-invariant by the number of square classes of the residue
class field (which for $\BR$ is $2$) similar arguments as in the
proof of Theorem 2 give as well an $\forall \exists$-formula.}
\begin{itemize} \item $k$ is $\BC (n = 1)$; \item $k$ is a
finite field or the function field of a real or complex curve $(n =
2)$; \item $k$ is a finite number field $(n=4)$; \item $k$ is the
function field of a complex variety $V ( n = u (k) \leq 2^d$ where
$d$ is the dimension of $V$).
\end{itemize}

There are, however, fields $k$ such that $k [[t]]$ is not $\forall
\exists$-definable in $k((t))$. Here is an example suggested by A.
Fehm: \\
Let $k^* = \bigcup_{n \geq 1} k_n$ with $k$ arbitrary and $k_n = k
((t_n)) \ldots (( t_1 ))$. \newline Clearly $k_n \subseteq k_{n+1}$.
It is also clear that $k^*$ is isomorphic to $K = k^* ((t))$ by
sending $t_1$ to $t$ and $t_{n+1}$ to $t_n$ for $n \geq 1$. The
pre-image $\fo$ of the henselian valuation ring $\CO = k^* [[t]]$ of
$K$ is again a henselian valuation ring of $k^*$. Note that the
restriction of $\CO$ to the subfield $k^*$ of $K$ is the trivial
valuation on $k^*$. \\

We now use the fact that $k^*$ is existentially closed in $K$ (see
Proposition 2' in \cite{E}). We then see that $\CO$ cannot be
$\forall \exists$-definable in $K$. Assume some $\forall
\exists$-formula $\varphi (x)$ would define $\CO$ in $K$. Then the
same formula would define $\fo$ in $k^*$. Now $t_1^{-1} \in \CO$
implies that $\varphi (t_1^{-1})$ hold in $K$. As $k^*$ is
existentially closed in $K$, we would also get $\varphi (t_1^{-1})$
in $k^*$. But then $t_1^{-1} \in
\fo$, a contradiction. \\

\section{$\exists$-definable henselian valuation rings}

Again we assume that $\Sigma$ is a first order axiom system for
henselian valued fields $(K, \CO)$. We consider two models $(K_1,
\CO_1), (K_2, \CO_2)$ of $\Sigma$ and assume $K_1 \subseteq K_2$ (as
fields). In order to get uniform $\exists$-definability for all
rings $\CO$ of the models $(K, \CO)$ of $\Sigma$, we have to show
that $\CO_1 \subseteq \CO_2$. Now we can no longer assume that $\CO
:= \CO_2 \cap K_1$ is a henselian valuation ring of $K_1$. Thus we
pass to the henselian closure $K^h$ of $K_1$ inside $K_2$ with
respect to $\CO_2$. The valuation $\CO_1$ being henselian on $K_1$
uniquely extends to a henselian valuation $\CO_{1}^{'}$ on $K^h$.
Thus we have now $\CO^h = \CO_2 \cap K^h$ and $\CO_{1}^{'}$ as
henselian valuations on $K^h$. In the next applications we shall fix
condition such that Fact 1 yields comparability of $\CO^h$ and
$\CO_{1}^{'}$. This clearly implies comparability of $\CO$ and
$\CO_1$. Thus it
remains to exclude $\CO \subsetneqq \CO_1$. \\

\begin{Thm} Let $\Sigma$ be a first order axiom system for henselian
valued fields $(K, \CO)$ such that the r.c.f. $\CO \diagup M = k$ is
finite or PAC and the fixed polynomial $f(x) \in \BZ [X]$ has no
zero in $k$. Then there is an $\exists$-formula $\varphi_f$ defining
uniformly the rings $\CO$ of models $(K, \CO)$ of $\Sigma$.
\end{Thm}

\textit{Proof:} If $\CO^h$ and $\CO_{1}^{'}$ would not be comparable
by Fact 1 both had a separably closed r.c.f. But then $f$ had a zero
in $\CO_2 \diagup M_2$. Thus we get comparability of $\CO$ and
$\CO_1$. We want to exclude $\CO \subsetneqq \CO_1$. In case it
holds, $\CO \diagup M_1$ is a proper valuation of $\CO_1 \diagup
M_1$. It then follows from Corollary 11.5.5 in \cite{F-J}, that the
r.c.f. $\CO \diagup M$ being the r.c.f. of $\CO \diagup M_1$ w.r.t.
$M \diagup M_1$ is separably closed. (Note that henselianity of $\CO
\diagup M_1$ is not needed.) But then again $f$ would have a zero in
$\CO_2
\diagup M_2$, a contradiction. \hfill{$\square$} \\

The result of Theorem 3 is due to A. Fehm. In \cite{F} he explicitly
constructs an $\exists$-formula $\varphi_f$. \\

\strut \\

\strut \\
\strut \\
\strut \\
\strut \\

Department of Mathematics and Statistics\newline University of
Konstanz, Konstanz, Germany
\newline alex.prestel@uni-konstanz.de


\begin{thebibliography}{quadmod}

\bibitem[A-K]{A-K}
W. Anscombe, J. Koenigsmann: \textit{An existential
$\emptyset$-definition of $\BF_q[[t]]$ in $\BF_q ((t))$}, ArXiv:
1306.6/60v1.

\bibitem[A1]{A1}
J. Ax: \textit{On the undecidability of power series fields}, Proc.
Amer. Math. Soc 16:846 (1965).

\bibitem[A2]{A2}
J. Ax: \textit{The elementary theory of finite fields}, Ann. of
Math. 88 (1968), 239-271.

\bibitem[C-D-L-M]{C-D-L-M}
R. Cluckers, J. Deraskhshan, E. Leenknegt, A. Macintyre:
\textit{Uniformly defining valuation rings in henselian valued
fields with finite or pseudo-finite residue fields}, Preprint.

\bibitem[E-K-M]{E-K-M}
R. Elman, N. Karponko, A. Merkurjev: \textit{The algebraic and
geoemtric theory of quadratic forms}, AMS, Coll Pub, 56 (2008).

\bibitem[E-P]{E-P}
A. Engler, A. Prestel: \textit{Valued Fields}, Springer Monographs
in Mathematics (2005).

\bibitem[E]{E} Y.L. Ershov: \textit{Relative regular
closedness and $\pi$-valuations}, Algebra and Logic, 31(6) (1992),
140-146.

\bibitem[F]{F} A. Fehm: \textit{Existential
$\emptyset$-definability of henselian valuation rings}, Preprint,
ArXiv 1307.1956v2.

\bibitem[F-J]{F-J}
M.D. Fried, M. Jarden: \textit{Field arithmetic}, 2nd Ed., Springer
Ergebnisse (2005).

\bibitem[P-D]{P-D}
A. Prestel, C.N. Delzell: \textit{Positive polynomials}, Springer
Monographs in Mathematics (2001).


\end{thebibliography}
\end{document}